\documentclass[12pt]{amsart}
\usepackage[utf8]{inputenc}

\usepackage{subfiles}
\usepackage{hyperref}
\usepackage{cleveref}
\usepackage{amsmath}
\usepackage{amssymb}
\usepackage{amsthm}
\usepackage{mathtools}
\usepackage{thmtools}
\usepackage{thm-restate}
\usepackage[margin=1in]{geometry}
\usepackage{cleveref}
\usepackage{ytableau}
\usepackage{comment}
\usepackage{soul}
\usepackage{tikz}
\usetikzlibrary{positioning}

\usepackage{graphicx}
\usepackage{caption}
\usepackage{float}

\usepackage[maxnames=5]{biblatex}
\renewbibmacro{in:}{}
\DeclareFieldFormat[article]{citetitle}{#1}
\DeclareFieldFormat[article]{title}{#1}

\numberwithin{equation}{section}

\theoremstyle{definition}
\newtheorem{theorem}{Theorem}[section]
\newtheorem*{theorem*}{Theorem}
\newtheorem{example}[theorem]{Example}
\newtheorem*{example*}{Example}
\newtheorem{lemma}[theorem]{Lemma}
\newtheorem*{lemma*}{Lemma}

\newtheorem*{corollary*}{Corollary}

\newtheorem*{definition*}{Definition}
\newtheorem{proposition}[theorem]{Proposition}
\newtheorem*{proposition*}{Proposition}
\newtheorem{remark}[theorem]{Remark}
\newtheorem*{remark*}{Remark}
\newtheorem{conjecture}[theorem]{Conjecture}

\usepackage{ mathrsfs }

\title{More results on stack-sorting for set partitions}

\author{Samanyu Ganesh}\address{\textsc{S. Ganesh}, The Westminster Schools, Atlanta, GA, 30327} \email{samanyuganesh@westminster.net}
\author{Lanxuan Xia}\address{\textsc{L. Xia}, St. Mark's School,
    Southborough, MA, 01772} \email{lanxuanxia@stmarksschool.org}
\author{Bole Ying}\address{\textsc{B. Ying}, Lower Merion High School, Ardmore, PA, 19003} \email{s017919@students.lmsd.org}
\begin{document}

\maketitle

\begin{abstract}
\label{abstract}
Let a \textit{sock} be an element of an ordered finite alphabet A and a sequence of these elements be a \textit{sock sequence}. In 2023, Xia introduced a deterministic version of Defant and Kravitz's stack-sorting map by defining the $\phi_{\sigma}$ and $\phi_{\overline{\sigma}}$ pattern-avoidance stack-sorting maps for sock sequences. Xia showed that the $\phi_{aba}$ map is the only one that eventually sorts all set partitions; in this paper, we prove deeper results regarding $\phi_{aba}$ and $\phi_{\overline{aba}}$ as a natural next step. We newly define two algorithms with time complexity $O(n^3)$ that determine if any given sock sequence is in the image of $\phi_{aba}$ or $\phi_{\overline{aba}}$ respectively. We also show that the maximum number of preimages that a sock sequence of length $n$ has grows at least exponentially under both the $\phi_{aba}$ and $\phi_{\overline{aba}}$ maps. Additionally, after defining fertility numbers (which were introduced by Defant in his study of stack-sorting for permutations) in the context of set partitions, we show that every positive integer is a fertility number under $\phi_{aba}$ and $\phi_{\overline{aba}}$. We prove that any multiple-pattern-avoiding stack for which a pattern order-isomorphic to $aba$ is not avoided will not sort all sock sequences. Finally, we show that there are exactly $2^{n-1}$ sock sequences of length $n$ that are 1-stack-sortable under $\phi_{aba, aab}$, where $n$ is any positive integer.                                                                                                                                                                                                                                                                                                                                                                                                                                                                                                                                                                                                                                                                                                                                                                                                                                                                           
\end{abstract}

\section{Introduction} 
\label{intro}
In 1968, Knuth \cite{knu2011} began the study of permutation patterns with a non-deterministic stack-sorting machine. A stack is a data structure that follows a first-in-last-out principle. Adding an element to the stack is referred to as a \textit{push} and removing an element from the stack is referred to as a \textit{pop}. Knuth \cite{knu2011} proved that a permutation $p$ could be sorted by the stack with the correct sequence of pushes and pops if and only if it avoids the $231$-pattern. Furthermore, Knuth \cite{knu2011} showed that the number of sortable permutations of length $n$ is given by the $n^{\text{th}}$ Catalan number $C_n=\frac{1}{n+1}{2n\choose n}$.

In 1990, West \cite{wes1990} proposed a deterministic version of Knuth’s stack-sorting machine. In West's stack-sorting algorithm, a push occurs when the next entry is smaller than the topmost element in the stack, and a pop occurs otherwise. In 2020, Cerbai, Claesson, and Ferrari \cite{cerclaefer} introduced pattern-avoiding stack-sorting algorithms, a generalization of West's algorithm that requires that a specific pattern is not formed by the stack's elements in relative order from top to bottom. Also, in 2020, Defant and Zheng \cite{defzhe2020} introduced consecutive-pattern-avoiding stack-sorting algorithms, in which the restricted pattern must not appear as a series of consecutive elements in the stack, again from top to bottom. 
 
More recently, in 2022, Defant and Kravitz \cite{defkrav2022} proposed an analog of Knuth's nondeterministic stack-sorting machine for set partitions. Then in 2023, Xia \cite{xia2023} introduced a deterministic version of Defant and Kravitz's stack-sorting map by defining the $\phi_{\sigma}$ and $\phi_{\overline{\sigma}}$ pattern-avoidance stack-sorting maps, where $\sigma$ denotes a pattern of socks to be avoided. For $\phi_{\sigma}$, a push occurs unless it would create a subsequence order-isomorphic to $\sigma$ (read from top to bottom in the stack); otherwise, a pop occurs. Similarly, for $\phi_{\overline{\sigma}}$, a push occurs unless it would create a subsequence order-isomorphic to $\sigma$ as a consecutive pattern in the stack, again read from top to bottom. Xia \cite{xia2023} proved that the $\phi_{aba}$ map is the only one that eventually sorts all set partitions. Naturally, we look further into $\phi_{aba}$ whilst also proving analogous results regarding $\phi_{\overline{aba}}$.

Formally, $\phi_{aba}$ dictates that at every step, a push is performed unless a pattern order-isomorphic to $aba$ would be created in the stack, in which case a pop is performed. For $\phi_{\overline{aba}}$, the condition is less strict, as the socks that comprise a potential $\overline{aba}$ pattern must appear consecutively in the stack.

Additionally, we newly define \textit{multiple-pattern-avoidance} stack-sorting algorithms for set partitions, which are analogous to Berlow's \cite{ber2021} multiple-pattern-avoidance stack-sorting algorithm for permutations in which the stack avoids multiple different patterns at the same time. Let $S=\{\sigma_1,\dots,\sigma_n\}$ be a collection of sock patterns. Then, $\phi_S$ denotes the stack-sorting algorithm that, at every step, requires a push unless pushing would create a subsequence that is order-isomorphic to any pattern in $S$ (when read from top to bottom in the stack) and otherwise requires a pop. It follows from Xia's findings that if $aba\in S$, $\phi_S$ eventually sorts all set partitions.

Xia \cite{xia2023} introduced the following definition of $k$-stack-sortability for set partitions: a set partition is $k$-stack-sortable if and only if it first becomes sorted after $k$ passes through the stack. She found a generating function that enumerates the number of 1-stack-sortable sequences of any given length that contain a fixed number of unique socks. In 2024, Choi, Gan, Li, and Zhu \cite{choi2024set} characterized the set partitions of minimal length and next minimal length that require the maximum number of passes through the $aba$-avoiding stack to be sorted. Later, Zhang \cite{Zhang2024} introduced minimally-sorted permutations as well as strengthening Berlow's \cite{ber2021} classification of periodic points. 

\begin{figure}
    \hspace{-1.25cm}\includegraphics[scale=0.4]{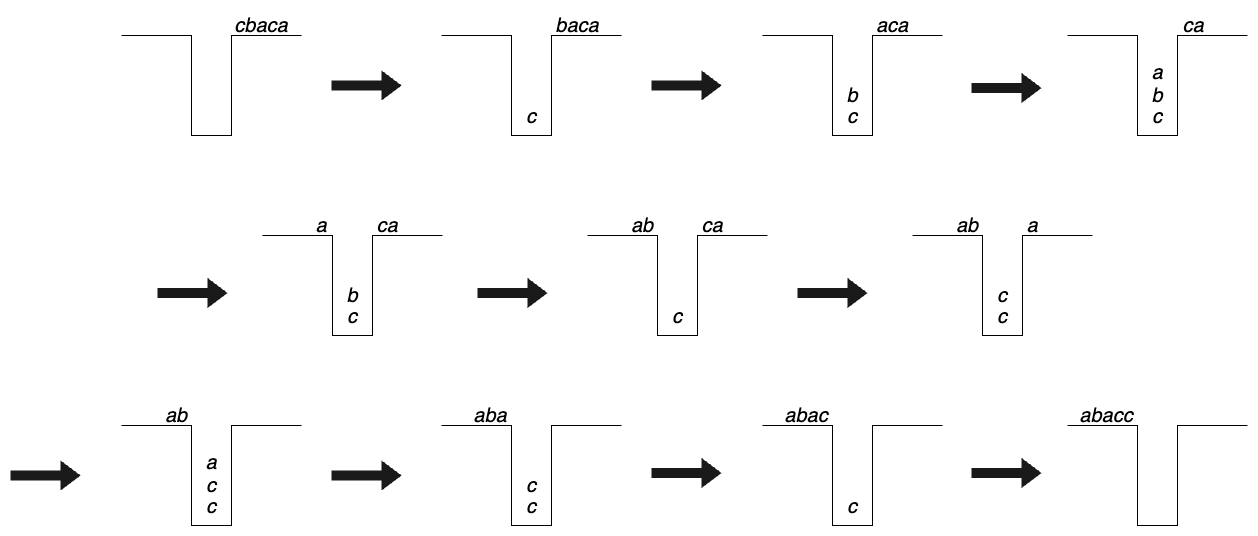}
    \begin{center}
        \caption{Applying $\phi_{aba}$ to the set partition $p = cbaca$}
        \label{fig:enter-label}
    \end{center}
\end{figure}

In \Cref{proofs1}, we define and prove the accuracy of two algorithms, both of which run in $O(n^3)$ time: the first determines if a given set partition is in the image of $\phi_{\overline{aba}}$, and the second performs the same function with regard to the $\phi_{aba}$. \Cref{proofs2} centers on the preimages and fertility numbers of both the $\phi_{\overline{aba}}$ and $\phi_{aba}$ maps. More specifically, we show that the maximum number of preimages possessed by sequences of length $n$ grows at least exponentially with $n$ for $\phi_{\overline{aba}}$ as well as for $\phi_{aba}$. After defining fertility numbers in the context of set partitions in \Cref{prelim} and again in 
\Cref{proofs2}, we show that every positive integer is a fertility number of both $\phi_{\overline{aba}}$ and $\phi_{aba}$ and contrast this result with an analogous statement proved in the context of stack-sorting for permutations. Finally, \Cref{proofs3} deals with multiple-pattern-avoidance algorithms: we generalize a result originally shown by Xia \cite{xia2023} for single-pattern-avoidance algorithms and derive an explicit formula for the number of set partitions of a given length that are 1-stack-sortable under $\phi_{aba, aab}$. Our final result launches a brief discussion regarding the multiple-pattern-avoidance algorithms that yield interesting sets of 1-stack-sortable sequences.

\section{Preliminaries}
\label{prelim}
We map every set partition of $[n]$ to an equivalence class of \textit{sock sequences}, where \textit{socks} are elements of an arbitrary finite alphabet $A$ that has an inherent order. Throughout the rest of this paper, we will be using the traditional $\{a,b,c,\dots\}$ or $\{a_1,a_2,\dots\}$ to represent distinct socks. For example, the set partition $\{\{1, 3\}, \{2, 5\}, \{4\}\}$ corresponds to the sock sequences $abacb$ or $a_1a_2a_1a_3a_2$, depending on the alphabet used. 

A sorted sock sequence is one where all occurrences of a particular element appear consecutively. For any sock sequence, we \textit{standardize} it by injectively renaming its elements such that the first occurrences of each element from left to right are assigned socks that follow the inherent order of the alphabet $A$. If the standardized sequences are the same, we call the original sequences \textit{equivalent}. To provide an example, the sequences $p_1=baacdcc$ and $p_2=wyyzxzz$ are equivalent up to standardization since both sequences become $p=abbcdcc$ after standardization using the alphabet $A=\{a,b,c,\ldots\}$. 

An \textit{equivalence class} of sock sequences is one where all sequences in the class are equivalent up to standardization. Additionally, if $M=\{a_1,\dots,a_n\}$ is a set (or multiset) of socks, denote $S(M)$ to be all possible sock sequences formed by elements of $M$, up to equivalence. For example, $a_1, a_1a_2, a_1a_1a_2,\dots\in S(M)$, but $a_1a_2$ and $a_2a_1$ cannot both be in $S(M)$ since they are equivalent by our definition. 

We can instead consider both the $\phi_{\overline{aba}}$ and $\phi_{aba}$ algorithms as maps. For example, $\phi_{aba}$ can be considered a map from $A$ to itself. Hence, we refer to the result of applying $\phi_{aba}$ to a sequence $p$ as the \textit{image} of $p$ under the $\phi_{aba}$ map. Moreover, if $p^*$ is the image of $p$ under a stack-sorting map, we call $p$ the \textit{preimage} of $p^*$ under the $\phi_{aba}$ map. Note that these maps are not necessarily bijective, surjective, or injective. $p^*$ can have multiple preimages, and not all sock sequences are even in the image of $\phi_{\overline{aba}}$ or $\phi_{aba}$. In the next section, we present numerous results regarding the images and preimages of both maps.

As an analog of Defant's \cite{def2018} definition of fertility numbers for permutations, we define a \textit{fertility number under a given map $\phi_\sigma$} in the context of sock sequences as a positive integer $k$ for which there exists a sock sequence that has $k$ preimages under $\phi_\sigma$. 

We cite a lemma and a proposition from Xia's paper \cite{xia2023}. This gives a general characterization of the image of $\phi_{aba}$.

\begin{lemma}[Xia {\cite[Lemma 3.1]{xia2023}}]  \label{Xialemma2.1}
    Let $p=x^{l_1}s_1x^{l_2}s_2\cdots x^{l_m}s_mx^{l_{m+1}}$ be such that $l_1,l_2,\cdots,l_m>0$, $l_{m+1}\geq0$, and $s_1,s_2,\cdots,s_m$ contain no occurrences of $x$, then 
    \begin{align*}
        \phi_{aba}(p) = \phi_{aba}(s_1) \phi_{aba}(s_2) \cdots \phi_{aba}(s_m) x^{l_1 + l_2 + \cdots + l_{m+1}}. 
    \end{align*}
\end{lemma}

The next proposition shows that $aba$ is the only pattern that eventually sorts all sequences in the case of single sock avoidance. This result helps us develop the more general case of multiple sock avoidance.

\begin{proposition}[Xia {\cite[Proposition 5.2]{xia2023}}] \label{Xiaproposition2.2}
    For $\sigma$ not of the form $a\cdots aba\cdots a$ and any multiset $M$ of socks where not all sequences in $S(M)$ are sorted, there exists a sequence $p\in S(M)$ that is not k-stack-sortable under $\phi_{\sigma}$ for any k.
\end{proposition}

To conclude this section, we establish a few new definitions and prove an analog of Lemma 2.1 for the $\phi_{\overline{aba}}$ map. For any sock sequence, we define the \textit{premier sandwiched sock} as the first sock in the sequence adjacent to two appearances of the same sock on both sides, with the additional condition that the sock differs from its neighboring socks. Consider the following algorithm: begin with a sock sequence $x$ and an empty set $Q$. Then, remove the premier sandwiched sock of $x$ and add it to $Q$. Then, remove the premier sandwiched sock of the resulting sequence and add that to $Q$. Continue until the sock sequence that remains (a subset of the elements of the original sequence) has no premier sandwiched sock; let the sequence that remains be denoted $x^\prime$ and the final set $Q$ be denoted the \textit{sandwich set} of $x$. Finally, let $\textnormal{rev}(x)$ denote the sock sequence obtained by reversing the socks in a sequence $x$.
\begin{lemma} \label{Lemma2.3}
    Let $Q=\{\eta_1, \eta_2,\ldots,\eta_k\}$ denote the sandwich set of a particular sock sequence $p$, where the $\eta_i$ are in order that they appear in $p$. Then
        $$\phi_{\overline{aba}}(p)=\eta_1 \eta_2\ldots\eta_k\textnormal{rev}(p^\prime).$$
\end{lemma}
\begin{proof}
     By definition, we push the socks of $p$ onto the stack until we reach $\eta_1$. We push $\eta_1$ onto the stack, but then we must pop $\eta_1$ since pushing the next sock onto the stack would create an $aba$-consecutive pattern.
    
    We then continue pushing elements of $p$ onto the stack until the next sock in the sandwich set is reached, whereupon we follow the same procedure outlined above (popping each $\eta_i$ and then continuing on). Once all of the socks $\eta_i$ have been popped, the remaining socks of the stack form the sequence $p^\prime$, from bottom to top. We pop the entire stack to finish, which appends $\textnormal{rev}(p^\prime)$ to the end of the output sequence.
\end{proof}



\section{Characterizing the images of $\phi_{\overline{aba}}$ and $\phi_{aba}$}
\label{proofs1}
In our first two theorems, we present two polynomial-time algorithms: the first determines whether a given set partition is in the image of $\phi_{\overline{aba}}$, and the second determines whether it is in the image of $\phi_{aba}$.
Note that passing all possible set partitions of length $n$ through the stack and comparing each output to the desired set partition runs in exponential time, as the number of set partitions of length $n$ is given by the Bell numbers (which grow exponentially fast).
Thus, our results significantly improve upon existing methods for the classification set partitions into the images of $\phi_{\overline{aba}}$ and $\phi_{aba}$.

As mentioned above, we first introduce the algorithm for the simpler $\phi_{\overline{aba}}$ map.

\begin{theorem} \label{cons_aba_image}
    There exists an algorithm that determines if any given set partition is in the image of $\phi_{\overline{aba}}$ and runs in $O(n^3)$ time.
\end{theorem}

\begin{proof}
    We first present the algorithm:
    \begin{enumerate}
    \item Iterate through the sequence, checking if everything to the right of the current sock is already sorted.
    \item If so, denote everything to the left of the current sock as the \textit{left-subsequence} and the rest of the sequence as the \textit{right-subsequence}. Consider the first sock in the left-subsequence. Iterate through the right-subsequence in reverse, looking for two consecutive socks that (a) are the same as each other and (b) are not the same as the first sock in the left-subsequence. Let these socks be denoted as a \textit{consecutive pair}.
    \item  If no such consecutive pair can be found, then go back to (1) and find the next index at which the overall sequence can be split such that the right-subsequence is already sorted.
    \item If such a pair is found, repeat (2) for the next sock in the left-subsequence, searching for two consecutive socks in the right-subsequence that (a) are the same as each other and (b) are not the same as the currently chosen sock in the left-subsequence. However, the search for this consecutive pair must be done by iterating through the right-subsequence in reverse (just as before), with the stipulation that it must begin where we left off previously (i.e. begin with the leftmost sock of the \textit{previous} consecutive pair that was found).
    \item If the end of the left-subsequence is reached (unique consecutive pairs are found for each of them, in reverse order within the right-subsequence), then the given sock sequence is in the image.
    \item If all possible splits of the given sequence have been considered and the end of the left-subsequence has never been reached, the given sock sequence is not in the image.
    \end{enumerate}

    We now show that this algorithm successfully classifies all set partitions, which is evident by \Cref{Lemma2.3}. The algorithm yields a positive result if and only if there exists a split of the sequence such that every sock in the left-subsequence has a corresponding consecutive pair in the right-subsequence, with these corresponding pairs appearing in order from left to right. Thus, a positive result from the algorithm implies we can construct the preimage of the sequence using \Cref{Lemma2.3}: this is done by reversing the right-subsequence and then placing the socks in the left-subsequence in between their corresponding pair in the right-subsequence.

    Next, we show that the algorithm runs in $O(n^3)$ time. Searching for a split such that the right-subsequence is already sorted runs in quadratic time, since iterating through each possible split of the sequence is obviously linear and the process of checking whether the right-subsequence is already sorted is also linear: to do the latter, we simply iterate through from left to right, keeping a list of the socks that have already been seen and ensuring that each sock encountered either has not already been seen or is the same as the one that directly precedes it. 
    
    Subsequently, the creation of a map between socks in the left-subsequence and consecutive pairs in the right-subsequence runs in $O(n)$ time. This process is linear due to the simple fact that, as we iterate through the left-subsequence and search for consecutive pairs in the right-subsequence, we only end up iterating (in reverse) through the right-subsequence once.

    Overall, since these two independent steps of the algorithm run quadratic and linear time, the overall algorithm is cubic in nature.
\end{proof}
The next theorem is significant in that it characterizes the image of $\phi_{aba}$ with the same time complexity as the algorithm outlined above that deals with $\phi_{\overline{aba}}$; \Cref{Xialemma2.1} is recursive in nature, so the existence of a non-recursive, polynomial-time algorithm that determines if a given sock sequence is in the image of $\phi_{aba}$ is surprising. 

\begin{theorem} \label{aba_image}
    There exists an algorithm that determines if any given set partition is in the image of $\phi_{aba}$ and runs in $O(n^3)$ time.
\end{theorem}
\begin{proof}
We begin by illustrating the algorithm's steps along with two example sock sequences: $\textit{bcbabccdd}$, which \textbf{is} in the image, and $\textit{bcbcbaabcccdd}$, which \textbf{is not} in the image.
\begin{enumerate}
        \item Iterate through the given sequence from left to right until an $aba$ pattern is found. Then, place a \textit{divider} right directly before the second appearance of the \textit{a} in the $aba$ pattern that was found. Repeat this process until no $aba$ patterns can be found to the right of the most-recently-placed divider. For the sake of visualization, these dividers will be denoted with the ``$\|$" symbol. 
        \item After this procedure, our examples would look as follows:
        \begin{center}
            $\textit{bcbabccdd} \Longrightarrow \textit{bc$\|$ba$\|$bccdd}$\\
            $\textit{bcbcbaabcccdd} \Longrightarrow \textit{bc$\|$bc$\|$baa$\|$bcccdd}$
        \end{center}
    \item Now, we establish a function $\gamma(i)$ that takes in a sock sequence $i$ as input and returns an integer. The algorithm consists of moving from left to right across the sequence and updating $\gamma(i)$. We claim that a sequence $i$ is in the image if and only if $\gamma(i)\geq0$.
    \item The dividers separate the sequence into what we refer to as \textit{blocks}. $\gamma(i)=0$ to start. There are three rules as to how $\gamma(i)$ changes and dividers are manipulated:
    \begin{enumerate}
        \item Every time a divider is encountered, $\gamma(i)$ decreases by one.
        \item As we move from left to right across the sequence, consider the strings of the same sock that appear consecutively. Every time one such run of the same sock is encountered, $\gamma(i)$ is updated. Consider one such run of socks that has length $l\geq1$ and appears in the $j$-th block of the sequence (counting from the left). Let the repeated sock be denoted $a$. Then, starting from the sock directly preceding $a$ and moving leftwards, consider each sock in the sequence until another $a$ is encountered. Let this first $a$ appear in the $m$-th block, counting from the left; if there are no appearances of $a$, then let $m=0$. Let $k=\min(l-1,j-m-1)$, then $\gamma(i)$ increases by $k$. We finish by manipulating dividers in accordance with the value of $k$: 
        \begin{itemize}
            \item If $k>0$, get rid of the dividers separating the $(j-k)$-th through the $j$-the block, effectively combining them into one new, larger block.
            \item If $k=0$, neither add nor remove dividers.
            \item If $k=-1$, place a divider directly before the current sock.
        \end{itemize}
    \end{enumerate}

    \item Ultimately, if $\gamma(i)\geq0$ upon the algorithm's completion, the sock sequence $i$ is in the image of the $\phi_{aba}$ map. Contrarily, a negative final value for $\gamma(i)$ indicates that the sock sequence $i$ is not in the image of the $\phi_{aba}$ map.

    \item The below table illustrates the value of $\gamma(i)$ for our two examples at each point where $\gamma(i)$ changes by a nonzero amount. The bolded characters illustrate where we are in the sequence at a particular step (note that if a sock that directly follows a divider is bolded, it means that one was just subtracted from $\gamma(i)$ due to the divider's presence).

    \begin{center}
    \begin{tabular}{ c | c | c | c }
    \hline
    Sequence & $\gamma(i)$ & Sequence & $\gamma(i)$ \\
    \hline
     \textit{\textbf{b}c$\|$ba$\|$bccdd} & 0 & \textit{\textbf{b}c$\|$bc$\|$baa$\|$bcccdd} & 0\\ 
     \textit{bc$\|$\textbf{b}a$\|$bccdd} & -1 & \textit{bc$\|$\textbf{b}c$\|$baa$\|$bcccdd} & -1\\ 
     \textit{bc$\|$ba$\|$\textbf{b}ccdd} & -2 & \textit{bc$\|$bc$\|$\textbf{b}aa$\|$bcccdd} & -2\\ 
     \textit{bc$\|$babc\textbf{c}dd} & -1 & \textit{bc$\|$bcba\textbf{a}$\|$bcccdd} & -1\\  
     \textit{bcbabccd\textbf{d}} & \boxed{0} & \textit{bc$\|$bcbaa$\|$\textbf{b}cccdd} & -2\\   
     & & \textit{bc$\|$bcbaa$\|$bcc\textbf{c}dd} & -2\\  
     & & \textit{bc$\|$bcbaabcccd\textbf{d}} & \boxed{-1}\\  
    \end{tabular}
    \end{center}
\end{enumerate}

We now show that this algorithm successfully classifies all set partitions, which is made easier by \Cref{Xialemma2.1}.
The initial assignment of dividers splits up the sequence into the least number of blocks such that each block is already sorted (i.e. does not contain any patterns that are order-isomorphic to $aba$). It is obvious from \Cref{Xialemma2.1} that any sorted sock sequence is in the image of $\phi_{aba}$. Essentially, the algorithm attempts to combine these initial blocks (by removing dividers) into longer and longer sequences that are also in the image of $\phi_{aba}$.

If $p$ blocks $r_1, r_2,\ldots,r_p$ appear in that order, then the definition of dividers implies that $r_1,r_2,\ldots,r_p\in \text{img}(\phi_{aba})$ but $r_1r_2\ldots r_{p-1} \notin \text{img}(\phi_{aba})$. The recursive nature of \Cref{Xialemma2.1} tells us that combining $r_1,r_2,\ldots, r_{p-1}$ whilst ensuring that the result is in the image of $\phi_{aba}$ must use $p-1$ consecutive copies of the same sock $x$, which we will refer to as the \textit{splitting block}. \Cref{Xialemma2.1} requires that the sock $x$ must not be found in $r_1,r_2,\ldots r_{p-1}$; note that it is impossible for $x$ to appear at the end of $r_{p-1}$ while $r_p$ begins with $x^{p-1}$ since dividers will never separate two copies of the same sock due to the manner in which they are assigned. Applying \Cref{Xialemma2.1} with this setup shows that $r_1r_2\ldots r_{p-1}x^{p-1}$ has $\Pi_{k=1}^pxr_k$ as its preimage and thus must be in the image of $\phi_{aba}$. So, we combine the blocks $r_1, r_2,\ldots,r_p$ (removing $p-1$ dividers in the process). However, when any $p$ blocks $r_1, r_2,\ldots,r_p$ are combined according to the algorithm, the combined sequence $r_1 r_2\ldots r_{p-1}r_p$ is not yet confirmed to be in $\text{img}(\phi_{aba})$ until the algorithm has passed completely through the end of $r_p$. Consider the sock $y$ appearing after the splitting block in $r_p$ that also appears in one of $r_1,r_2,\ldots,r_{p-1}$. By \Cref{Xialemma2.1}, $r_1 r_2\ldots r_{p-1}r_p$ cannot be in $\text{img}(\phi_{aba})$; we handle this scenario with the $k=-1$ case of step 5b: a divider is added directly before the sock $y$ and now the block directly preceding this divider (which contains the original $r_1,r_2,\ldots,r_{p-1}$, and part of the original $r_p$) is in $\text{img}(\phi_{aba})$.

Thus, once the algorithm completes, all of the blocks that constitute the sequence must be in $\text{img}(\phi_{aba})$. So, if no dividers remain ($\gamma(i)\geq0$), then we have proven that the entire sequence is in $\text{img}(\phi_{aba})$. 

To prove that $\gamma(i)<0$ implies that the sequence is not in $\text{img}(\phi_{aba})$, we notice that the algorithm maximizes the number of dividers removed by considering every possible splitting block as it passes through the sequence. The theoretical maximum for the number of dividers removed would simply be the sum of the lengths (minus one) of each consecutive string of socks (i.e. $(1-1)+(3-1)+(2-1)=3$ for $abbbcc$). The algorithm cannot attain this maximum when the socks in the splitting block also appear in a previous block, which is handled in step 5: we simply remove the maximum number of dividers such that the combined block does not contain non-consecutive copies of the sock in the splitting block.

Finally, we show that the time complexity of the algorithm is $O(n^3)$. The initial assignment of dividers is a linear process: we keep track of the already-seen socks as we move through the sequence, and every time a pattern isomorphic to $aba$ is discovered, a divider is placed and we begin iterating immediately after the divider.

The rest of the algorithm consists of examining consecutive strings of the same sock $a$. Finding these strings is a linear process (we keep track of the most recently seen sock as we iterate through the sequence). Once each string is found, we iterate in reverse from the beginning of the string to the beginning of the entire sequence, looking for another appearance of $a$ in a preceding block that will eventually determine the value of $k$ (to compute this, we also keep track of where dividers are placed at any given time in order to count blocks easily)—this is also a linear process.

Thus, the procedures described in steps 4 and 5 run in quadratic time. Combining this with the fact that the initial assignment of dividers is a linear process, the entire algorithm therefore runs in cubic time.
\end{proof}

\section{On the preimages and fertility numbers of $\phi_{\overline{aba}}$ and $\phi_{aba}$}
\label{proofs2}
This section first showcases a pair of results that pertain to the growth in the number of preimages of $\phi_{\overline{aba}}$ and $\phi_{aba}$ as the sequence length increases. Subsequently, we present two theorems on how every integer is a fertility number under $\phi_{\overline{aba}}$ and $\phi_{aba}$ respectively. 

We begin with the following result concerning preimages under the $\phi_{\overline{aba}}$ map.
\begin{theorem}\label{pconsimg_exp}
    The maximum number of preimages for sequences of length $n$ grows at least exponentially with $n$ for $\phi_{\overline{aba}}$.
\end{theorem}
\begin{proof}
     Consider the sock sequence $p^* = a_1a_2\cdots a_n a_{n+1}^{k}$ with length $l$. We construct one of its preimages $p$, which clearly must be a sock sequence of length $k+n$.
    The first sock of $p$ must be $a_{n+1}$. Now consider the $k+n-1$ remaining spots of $p$. We choose $k-1$ of these spots for the remaining copies of $a_{n+1}$. Since there is only one copy of each of $a_1,a_2,\cdots,a_{n}$and $(k+n-1)-(k-1)=n$ spots remaining in $p$, the rest of the construction is forced by \Cref{cons_aba_image}.
    To see why, let $j\leq n$ denote the number of unoccupied spots in $p$ that occur before the final copy of $a_{n+1}$. By \Cref{cons_aba_image}, these spots must be filled by $a_1,a_2,\cdots,a_{j}$ in that order from left to right. Then, \Cref{cons_aba_image} also dictates that $p$ must end with the string $a_na_{n-1}a_{n-2}\ldots a_{j+1}$ (as there must necessarily exist $n-j$ unoccupied spots in $p$ that occur consecutively after the last copy of $a_{n-1}$). 
    Thus, there are exactly $\binom{k+n-1}{k-1}$ preimages of $p^*$. Then, setting $n=k=\frac{l}{2}$ results in $\binom{l-1}{\frac{l}{2}-1}$ preimages. Through Stirling's approximation, $\binom{l-1}{\frac{l}{2}-1} \approx \frac{2^n}{\sqrt{2\pi n}}$ so the maximum number of preimages of a sock sequence of length $n$ under the $\phi_{\overline{aba}}$ map grows at least exponentially with $n$. 
\end{proof}
We use a similar technique to prove the same result for $\phi_{aba}$.
\begin{theorem} \label{pimg_exp}
    The maximum number of preimages for sequences of length $n$ grows at least exponentially with $n$ for $\phi_{aba}$. 
\end{theorem}
\begin{proof}
    Consider the sock sequence $p^*=a_1a_2\cdots a_n a_{n+1}^k$ with length $l$. We construct one of its preimages $p$, which clearly must be a sock sequence of length $k+n$.
    Firstly, note that $p$ must start with $a_{n+1}$ by \Cref{Xialemma2.1}. Now consider the $k+n-1$ remaining spots of $p$. We choose $k-1$ of these spots for the remaining copies of $a_{n+1}$. Since there is only one copy of each of $a_1,a_2,\cdots,a_{n}$and $(k+n-1)-(k-1)=n$ spots remaining in $p$, the rest of the construction is forced by \Cref{Xialemma2.1}: the non-occupied spots must be filled by $a_n,a_{n-1},\cdots,a_{1}$ in that order from left to right. Thus, there are exactly $\binom{k+n-1}{k-1}$ preimages of $p^*$. Then, setting $n=k=\frac{l}{2}$ results in $\binom{l-1}{\frac{l}{2}-1}$ preimages. Again by Stirling's approximation, $\binom{l-1}{\frac{l}{2}-1} \approx \frac{2^n}{\sqrt{2\pi n}}$ so the maximum number of preimages of a sock sequence of length $n$ under the $\phi_{aba}$ map also grows at least exponentially with $n$. 
\end{proof}

In 2018, Defant \cite{def2018} introduced fertility numbers in the context of stack-sorting for permutations. We call a positive integer $k$ a \textit{fertility number} $\phi_{aba}$ if and only if there exists a set partition with exactly $k$ preimages. We show that every positive integer is a fertility number under both  $\phi_{aba}$ and $\phi_{\overline{aba}}$; additionally, for any positive integer $f$, we prove the existence set partitions of $[n]$ with $f$ preimages for arbitrarily large $n$.

\begin{theorem} \label{fert_aba_cons}
    Every positive integer is a fertility number under $\phi_{\overline{aba}}$. Furthermore, for all positive integer values of $f$, there exists a set partition of $[n]$ with exactly $f$ preimages under $\phi_{\overline{aba}}$, provided $n$ is sufficiently large ($n\geq f+1$).
\end{theorem}

\begin{proof}
    Specifically, we show that $q=a_1{a_2}^ma_3a_4\ldots a_{n-m+1}$, a sock sequence of length $n$, has $m$ preimages, where $m,n\in \mathbb{Z}$, $m\leq n-1$, and $a_k$ and $a_j$ denote distinct socks if and only if $k\neq j$.
    
    To see this, consider all sock sequences that have the same distribution of socks (up to equivalence) as the final sock sequence $a_1{a_2}^ma_3a_4\ldots a_{n-m+1}$. Clearly, inputs without this property don't map to the desired sock sequence.
    
    It is impossible for $a_2$ to be in the sandwich set of $q$, since there is only one copy of each of the socks $a_1, a_3, a_4,\ldots a_{n-m+1}$. Then, by \Cref{Lemma2.3}, any preimage of $q$ must have either an empty sandwich set or $a_1$ as the only member of its sandwich set. Essentially, any preimage of $q$ must begin with $a_{n-f+1}a_{n-f}\ldots a_4a_3a_2$. 
    
    Then, the remainder of the input sequence must be some rearrangement of $a_1$ and $m-1$ copies of $a_2$. Thus, the sock sequence $q=a_1{a_2}^ma_3a_4\ldots a_{n-m+1}$ has $m$ preimages.
\end{proof}

We finish off this section by showing the equivalent result for $\phi_{aba}$.

\begin{theorem} \label{fert_aba}
    Every positive integer is a fertility number under $\phi_{aba}$. Furthermore, for all positive integer values of $f$, there exists a set partition of $[n]$ with exactly $f$ preimages under $\phi_{aba}$, provided $n$ is sufficiently large ($n\geq f+1$).
\end{theorem}

\begin{proof}
    Specifically, we show that $q=a_1a_2a_3\ldots a_{m-1}a_{m}a_ma_{m+1}a_{m+2}\ldots a_{n-2}a_{n-1}$, a sock sequence sock sequence of length $n$, has $m$ pre-images, where $m,n\in \mathbb{Z}$, $m\leq n-1$, and $a_k$ and $a_j$ denote distinct socks if and only if $k\neq j$.

    To see this, consider all sock sequences that have the same distribution of socks (up to equivalence) as the desired sock sequence $a_1a_2a_3\ldots a_{m-1}a_{m}a_ma_{m+1}a_{m+2}\ldots a_{n-2}a_{n-1}$. Clearly, inputs without this property don't map to the desired sock sequence and thus are not preimages of $p$.

    By \Cref{Xialemma2.1}, the string of socks $a_{n-1}\ldots a_{m+1}a_m$ must appear at the beginning of any preimage of $q$, in that order. Then, let $k$ equal \textit{one more than} the number of socks strictly between the first and second appearances of $a_m$; clearly, $k\leq m$. \Cref{Xialemma2.1} dictates that the string of socks $a_{k-1}a_{k-2}\ldots a_2a_1$ must appear, in that order, between the two appearances of $a_m$. Moreover, by \Cref{Xialemma2.1}, any preimage of $q$ must conclude with the string of socks $a_{m-1}a_{m-2}\ldots a_{k+1}a_k$.
    There are thus $m$ spots for the second appearance of $a_m$, as $k$ ranges from 1 to $m$ inclusive. Thus, since each value of $k$ corresponds to exactly one preimage, the sock sequence $q=a_1a_2a_3\ldots a_{m-1}a_{m}a_ma_{m+1}a_{m+2}\ldots a_{n-2}a_{n-1}$ has $m$ preimages.
\end{proof}

\begin{remark}
\Cref{fert_aba_cons} and \Cref{fert_aba} notably differ from past results regarding fertility numbers under West's stack-sorting map for permutations. Defant and Zheng \cite{def2018} showed that, for permutations, the set of fertility numbers contains every nonnegative integer congruent to 3 modulo 4 and that the lower asymptotic density of the set is at least $\approx0.7618$.
Defant and Zheng \cite{def2018} also conjectured that the smallest fertility number congruent to 3 modulo 4 is 95. 
\end{remark}

\section{Results regarding multiple-pattern-avoidance stack-sorting maps}
\label{proofs3}
This section serves as a follow-up to Xia's \cite{xia2023} discussion of multiple-pattern-avoidance stack-sorting maps. Specifically, we show that if the $aba$ pattern is not one of the patterns avoided by a multiple-pattern-avoiding stack, then the stack does not sort all set partitions in finite time. Additionally, we prove an analogous result to {\cite[Theorem 1.1]{xia2023}} for the sorting map $\phi_{\{aba, aab\}}$. 

We first prove a generalization of Xia's {\cite[Proposition 5.2]{xia2023}} for multiple-pattern-avoidance stack-sorting algorithms.

\begin{theorem} \label{multpat}
    Let $G=\{\sigma_1,\dots,\sigma_n\}$ be a collection of sock patterns such that none of the $\sigma_i$ is of the form $a\dots aba\dots a$ and our stack sorting algorithm avoids all elements of $G$ simultaneously. Then let, $M$ be a set of socks such that $S(M)$ is not all sorted. There exists an unsorted element $p$ of $S(M)$ such that $p$ can never be sorted by the sorting algorithm corresponding to $G$.
\end{theorem}

\begin{proof}
    We consider two cases.\\
    \underline{Case 1:} Either all elements of $G$ contain a pattern order-isomorphic to $abba$ or $abca$, or all do not. By applying \Cref{Xiaproposition2.2}, we can evidently construct a sequence $p$ that is never sorted by the stack.\\
    \underline{Case 2:} $G$ has elements that contain patterns order-isomorphic $abba$ or $abca$ as well as elements that contain patterns order-isomorphic to $abab$, $abac$, or $caba$. In this case, consider the sequence $p=a_1a_2a_1a_3a_1a_4\cdots a_1a_ma_1$, where there is a new distinct sock between every pair of $a_1$s. As this sequence passes through the stack, it becomes $p^*=a_1a_3a_1a_4a\cdots a_1a_ma_1a_2a_1$, which is equivalent to the original sequence $p$, so $p$ will never be sorted.  
\end{proof}

Next, we focus on a specific multiple-pattern-avoidance stack-sorting algorithm, $\phi_{\{aba,aab\}}$, and replicate Xia's {\cite[Theorem 1.1]{xia2023}} and {\cite[Theorem 1.2]{xia2023}} for this map. 

In general, (multiple-)pattern-avoidance stack-sorting algorithms that do not avoid a pattern order-isomorphic to $aba$ have uninteresting collections of $1$-stack-sortable sequences of length $n$. In fact, $\phi_{\{aba, aab\}}$ is one of the only such maps for which there exist $1$-stack-sortable sequences of a given length $n$, and thus we characterize and count these sequences in \Cref{gnrfct}.

\begin{theorem} \label{gnrfct}
    Let $s(n)$ denote the number of sock patterns of length $n$ that are $1$-stack-sortable under $\phi_{\{aba,aab\}}$. Then $s(n)=2^{n-1}$, up to the equivalence of sock patterns.
\end{theorem}

\begin{proof}
     Note that when passing the sock sequence $p=x^{l_1}s_1x^{l_2}s_2\dots x^{l_m}s_mx^{l_{m+1}}$ through $\phi_{\{aba,aab\}}$, the $aab$ avoidance will cause all but $1$ $x$ in $x^{l_1}$ to exit the stack before any of $s_1$ enters the stack. Hence, if $m>1$ or $l_1>1$, occurrences of $x$ will be separated, hence leaving the sequence unsorted. Therefore, the only $1$-stack-sortable sequences under $\phi_{\{aba,aab\}}$ take the form $p=xs_1x^l$ where $l\geq0$ and $s_1$ is $1$-stack sortable. To construct such $p$, note there is one unique way to construct a $1$-stack-sortable sequence of length $n$ from a $1$-stack-sortable sequence of length $m<n$. Also note that we artificially define length $0$ $1$-stack sortable sequences to be a thing, and naturally $s(0)=1$ because the empty sequence is $1$-stack-sortable. Suppose $p(n)=\sum s(n)x^n$, then every $s(n)$ is the sum of all preceding $s(m)$'s (including $m=0$), the fitting sequence is $p(n)=2^{n-1}x^n$ since $s(0)=s(1)=1$. 
\end{proof}

\begin{example}
    $s(2)=s(0)+s(1)=1+1$, and the corresponding $1$-stack-sortable sequences of length $2$ are $ab$ and $aa$. Here, $ab$ is uniquely determined by the only $1$-stack-sortable sequence of length $1$, $a$, and $aa$ is uniquely determined by the trivial $1$-stack-sortable sequence of length $0$.
\end{example}

Next we present a refinement of the theorem above like the one stated in {\cite[Theorem 1.2]{xia2023}}. Again, instead of having a generating function like the case for counting the number of $1$-stack-sortable sequences of length $n$ and having $r$ distinct socks under $\phi_{aba}$, the case for $\phi_{\{aba,aab\}}$ is much simpler. 

\begin{theorem} \label{rfm}
    Let $s(n,r)$ denote the number of sock patterns of length $n$ that are $1$-stack-sortable under $\phi_{\{aba,aab\}}$ and have exactly $r$ distinct socks. Then $s(n,r)={n\choose r-1}$, up to the equivalence of sock patterns.
\end{theorem}

\begin{proof}
    As mentioned in the proof for \Cref{gnrfct}, there is one unique way to construct a $1$-stack-sortable sequence $p_n$ of length $n$ from a $1$-stack-sortable sequence $p_m$ of length $m<n$, that is, to have a sock $s$ that is distinct from every sock in the $p_m$, and append $1$ occurrence of $s$ in front of $p_m$ and $n-m-1$ occurrences of $s$ behind $p_m$. Hence, one can construct an unique sequence counted in $s(n,r)$ from an unique sequence counted in $s(m,r-1)$ where $r-1\leq m\leq n-1$, so $s(n,r)=\sum_{i=r-1}^{n-1} s(i,r-1)$. Note that since $\sum_{i=r-1}^{n-2} s(i,r-1)=s(n-1,r)$, $s(n,r)=s(n-1,r-1)+s(n-1,r)$. Thus, $s(n,r)={n\choose r-1}$, which makes sense as $s(n)=2^{n-1}$ and $\sum_{i=0}^{n-1} {n-1 \choose i}=2^{n-1}$. (Note, the $s(n,r)$'s corresponds to the entries in the $n-1$th row of the Pascal's triangle).
\end{proof}

\begin{remark}
    As mentioned in the lead-up to \Cref{gnrfct}, $1$-stack-sortability is rather uninteresting to study for multiple-pattern-avoidance stack-sorting algorithms where a pattern order-isomorphic to $aba$ is not avoided. 
    To see why, consider an arbitrary sock sequence as $p=x^{l_1}s_1x^{l_2}s_2\dots x^{l_m}s_mx^{l_{m+1}}$; without avoiding the $aba$ pattern, the copies of $x$ will never occur together in a consecutive block after just $1$ sort. Furthermore, $\phi_{aba}$ is the only (multiple-)pattern-avoidance stack-sorting algorithm that has a rather complicated collection of $1$-stack sortable sequences. The set of $1$-stack-sortable sequences under any other multiple-pattern-avoidance stack-sorting map is significantly limited, as the restrictions imposed by the map are either too loose or too tight (the latter of which is true for $\phi_{\{aba,aab\}}$, the subject of Theorems \ref{gnrfct} and \ref{rfm}).
\end{remark}


\section{Future Directions}
\label{futuredirections}
We propose several directions of future study that are motivated by computational results.

First, we present two conjectures concerning the sequences of length $n$ with the maximum number of preimages under $\phi_{\overline{aba}}$ and $\phi_{aba}$.

\begin{conjecture}
   The maximum number of preimages under $\phi_{\overline{aba}}$ that a set partition of length $n$ can have is given by $\frac{2x}{1-2x+\sqrt{1-4x^3}}$. \textcolor{red}{}
\end{conjecture}

\begin{conjecture}
    The set partition of length $\frac{n(n+1)}{2}$ with the most number of preimages under $\phi_{aba}$ is given by $a_{1}^1 a_{2}^{2} \ldots a_{n}^n$. 
\end{conjecture}

We also suggest an upper bound on the size of the image of $\phi_{\overline{aba}}$ for a fixed sequence length.

\begin{conjecture}
    The number of set partitions of length $n$ in the image of $\phi_{\overline{aba}}$ is bounded above by $2B(n-1)$, where $B(n)$ is the $n^{\text{th}}$ Bell number. 
\end{conjecture}

More generally, it would be helpful to derive an explicit formula or generating function that gives the number of preimages of any sorted sock sequence under both $\phi_{\overline{aba}}$ and $\phi_{aba}$.

Additionally, our characterization of the image for $\phi_{aba}$ is an initial first step to enumerating the image, which will help answer Xia's Question $2$ about the number of $k$-stack-sortable sock patterns under $\phi_{aba}$ for $k>1$. Hence, it is natural to ask for an actual enumeration of the size of the image of $\phi_{aba}$.  

\section*{Acknowledgements}
The authors thank Yunseo Choi for suggesting the problem and providing helpful feedback. 

\nocite{*}
\printbibliography

\end{document}